\documentclass[12pt]{amsart}

\usepackage{amssymb}
\usepackage{amsmath}
\evensidemargin\oddsidemargin

\begin{document}

\baselineskip=18pt
\setcounter{page}{1}

\newtheorem{theorem}{Theorem}
\newtheorem{lemma}[theorem]{Lemma}
\newtheorem{proposition}[theorem]{Proposition}
\newtheorem{remark}[theorem]{Remark}
\newtheorem{remarks}[theorem]{Remarks}
\newtheorem{THA}{Theorem A\!\!}
\newtheorem{THB}{Theorem B\!\!}
\newtheorem{CORO}{Corollary\!\!}
\newtheorem{CONJ}{Conjecture\!\!}

\renewcommand{\theTHA}{}
\renewcommand{\theTHB}{}
\renewcommand{\theCORO}{}
\renewcommand{\theCONJ}{}

\newcommand{\eqnsection}{
\renewcommand{\theequation}{\thesection.\arabic{equation}}
    \makeatletter
    \csname  @addtoreset\endcsname{equation}{section}
    \makeatother}
\eqnsection


\def\e{{\mathbb E}}
\def\FF{\mathcal{F}}
\def\LL{\mathcal{L}}
\def\k{{\mathcal K}}
\def\p{{\mathbb P}}
\def\ptil{{\widehat\p}}
\def\r{{\mathbb R}}
\def\z{{\mathbb Z}}
\def\SS{\mathcal{S}}
\def\lacc{\left\{}
\def\lcr{\left[}
\def\lpa{\left(}
\def\lva{\left|}
\def\racc{\right\}}
\def\rcr{\right]}
\def\rpa{\right)}
\def\rva{\right|}
\def\eps{\varepsilon}
\def\Un{{\bf 1}}
\def\ee{\mathrm{e}}
\def\d{\, \mathrm{d}}
\def\qed{\hfill$\square$}
\def\elaw{\stackrel{d}{=}}
%
%

\title[Lower tails of homogeneous functionals]
      {The lower tail problem for homogeneous functionals of stable processes with no negative jumps}

\author[Thomas Simon]{Thomas Simon}

\address{Equipe d'Analyse et Probabilit\'es, Universit\'e d'Evry-Val d'Essonne, Boulevard Fran\c{c}ois Mitterrand, F-91025 Evry Cedex. {\em E-mail address}: {\tt tsimon@univ-evry.fr}}

\keywords{Fluctuating additive functional, lower tail probability, self-similar process, stable process, Wiener-Hopf factorization in two dimensions.}

\subjclass[2000]{60F99, 60G18, 60G52, 60J55}

\begin{abstract} Let $Z$ be a strictly $\alpha$-stable real L\'evy process
  $(\alpha\in (1,2])$ and $X$ be a fluctuating $\beta$-homogeneous additive
  functional of $Z$. We investigate the asymptotics of the first
  passage-time of $X$ above 1, and give a general upper bound. When $Z$
  has no negative jumps, we prove that this bound is optimal and does
  not depend on the homogeneity parameter $\beta$. This extends a
  result of Y. Isozaki \cite{Iso} and solves partially a conjecture of Z. Shi \cite{Shi}. 
\end{abstract}

\maketitle

\section{Introduction}

Let $X$ be a real process starting from 0 and $T = \inf\{ t> 0, \, X_t
> 1\}$ be its first passage time above level one. Studying the law of $T$
is a classical and important question in probability, with many
applications especially in finance and insurance. Since it is often
difficult to compute the exact distribution of $T,$ even when $X$ is
continuous, one is sometimes interested in the asymptotic behaviour of the ruin probability: 
\begin{equation}
\label{prob1}
\p\lcr T > t\rcr, \qquad t \to\infty.
\end{equation}
In a self-similar framework, this problem is equivalent to the estimate of the so-called lower tails of $X$:
\begin{equation}
\label{prob2}
\p\lcr \sup_{t\in [0, 1]}X_t < \eps  \rcr, \qquad \eps \to 0.
\end{equation}
The latter is a less classical question, but striking connections have
been made over the years between (\ref{prob2}) and  subjects as
different as e.g. increase points for L\'evy processes - see section VI.5 in
\cite{Ber}, Hausdorff dimension of regular points for the inviscid
Burgers equation with random initial impulse \cite{Si1, MK}, capture
times of Brownian pursuits and zeroes of random polynomials - see \cite{LS, Sh} and also the references therein for other relevant issues. 

Notice that despite formal resemblance, the estimate (\ref{prob2}) has less to do with the so-called small deviation (or small ball) problem for $X$, which deals with the asymptotics $$\p\lcr\sup_{t\in [0,1]} |X_t| < \eps\rcr, \qquad 
    \eps\to 0.$$ 
Indeed the latter is usually studied, in a
large deviation spirit, under the logarithmic scale, whereas in
(\ref{prob2}) the investigated speed of convergence speed is most of the time polynomial: one expects
a behaviour of the type
\begin{equation}
\label{prob3}
\p\lcr T > t\rcr\; =\; t^{-\theta +o(1)}, \qquad t \to\infty
\end{equation}
for some finite positive constant $\theta,$ whose exact value is actually relevant in several problems of statistical physics - see \cite{BDM, BG} and the numerous references therein. In this physical literature, the exponent $\theta$ is often mentioned as a {\em persistence} or {\em survival} exponent because in general, the probability in (\ref{prob3}) behaves modulo some $o(1)$ term like the probability that $X$ starting from some $x > 0$ remains positive up to time $t\to +\infty.$ This expected polynomial asymptotic forces to work under the natural scale, making the problem significantly more delicate than usual small ball probability estimates and actually very few explicit exponents are known, especially in a Non-Markovian framework. So far, up to some generalisation of integrated Brownian motion which will be detailed soon afterwards, the problem of computing $\theta$ seems to have been solved only for the three following self-similar processes: 

\vspace{2mm}

\begin{center}
\begin{tabular}{|c|c|c|}
\hline
$X$ & $\theta$ & Ref. \\
\hline
\hline
Integrated Brownian motion & 1/4 & \cite{MacK, Go} \\
\hline
L\'evy $\alpha-$stable process & $\rho$ & \cite{Bi} \\
\hline
Fractional Brownian motion & $1-H$ & \cite{Mo} \\
\hline
\end{tabular}
\end{center}

\vspace{2mm}

In the above, the parameter $\rho$ stands for the positivity parameter of the L\'evy $\alpha-$stable process, which will be recalled henceforth, whereas the parameter $H$ is the usual Hurst parameter of the fractional Brownian motion. Recently, two conjectures were stated on the explicit expression of $\theta$ for integrated processes, generalizing the value $1/4$ for integrated Brownian motion:

\vspace{2mm}

\begin{center}
\begin{tabular}{|c|c|c|}
\hline
$X$ & $\theta$ & Ref. \\
\hline
\hline
Integrated L\'evy $\alpha-$stable process & $(\alpha-1)_+/2\alpha$ & \cite{Shi} \\
\hline
Integrated fractional Brownian motion & $H(1-H)$ & \cite{MK} \\
\hline
\end{tabular}
\end{center}

\vspace{2mm}

In this paper we will be concerned with the first conjecture. More precisely we will consider the slightly more general situation where $X$ is a fluctuating $\beta$-homogeneous additive functional of a strictly $\alpha$-stable L\'evy process $Z$:
$$X_t\; :=\;A^{(\beta)}_t\; =\; \int_0^t \vert Z_s\vert^\beta {\rm sgn} (Z_s) \, \d s, \quad t\ge 0,$$ 
for $\alpha \in (1,2]$ and $\beta > -(\alpha +1)/2.$ More details and the reason why the above signed integral always makes sense will be given in the
next section. The process $\{ A^{(\beta)}_t, \; t\ge 0\}$ is $(1 +
\beta/\alpha)$-self-similar, but is not Markovian and has no 
stationary increments. In the 
case $\beta =1$, it is the area process associated with $Z$, in other words the integrated L\'evy $\alpha$-stable process:
$$A^{(1)}_t\; =\; \int_0^t Z_s \, \d s, \quad t\ge 0,$$
and is also an $\alpha$-stable process in the general sense of
\cite{ST}. However, in the case $\beta \neq 1,$ it is no more a stable
process. When $\beta =-1,$ the process $A^{(\beta)}$ is up to a
multiplicative constant the Hilbert transform of $\{ L(t,x), x\in\r\}$ the local time
process of $Z$ - see section V.2 in
\cite{Ber}. When $\beta < -1,$ it can be viewed as a  fractional
derivative of $\{ L(t,x), x\in\r\}$ - see section 2 in \cite{FG}. Our
first result is the following:

\begin{THA} With the above notations, there exists a finite positive constant $\k$ such that
\begin{equation}
\label{Main1}
\p[ T> t]\; \le\; \k \, t^{-(\alpha -1)/2\alpha}, \qquad t\to +\infty.
\end{equation}
Besides, the constant $\k$ is explicit when $\alpha = 2$ or $\beta = -1$.
\end{THA}

In the Brownian case $\alpha = 2,$ this result had been obtained previously by Isozaki \cite{Iso} for $\beta >0,$ with the help of a Wiener-Hopf factorization in two dimensions extending an identity of Spitzer-Rogozin in discrete time, and an asymptotic analysis of the Wiener-Hopf factors. We remark that these analytical arguments can be carried out to the Non-Gaussian case without much difficulty. Our main novelty consists in dealing with {\em non necessarily
symmetric} stable processes. 

In \cite{Iso} it is also proved that when $\alpha =2$, the exponent $(\alpha-1)/2\alpha = 1/4$ is optimal, in the sense that (\ref{Main1}) holds as well with an inequality in the other direction and a positive constant $\k' < \k.$ With the help of excursion theory, this result in the case $\alpha =2$ has been then significantly refined in \cite{IK}, where it is proved that
$$\p[ T> t]\; \sim\; \k'' \, t^{-1/4}, \qquad t\to +\infty$$
for some explicit constant $\k'' \in (\k',\k).$ At a less precise level, finding the right exponent $\theta$ in (\ref{prob3}) remains a tantalizing question in the Non-Gaussian case, because of the
jumps of $Z$ which make analytical arguments untractable. Results in the direction of Shi's conjecture had been given in \cite{DSS}, with unfortunately a serious gap in the proof. In this paper, we aim at tackling this problem in the spectrally positive situation:

\begin{THB} Suppose that $Z$ has no negative jumps. Then with the above notations, there exists a positive constant $k$ such that
$$\p[ T> t]\; \ge\; k \, t^{-(\alpha -1)/2\alpha} (\log t)^{-1/2\alpha}, \qquad t\to +\infty.$$
\end{THB}
Putting the two theorems together entails that $\theta = (\alpha-1)/2\alpha$ for the processes $A^{(\beta)}$ as soon as $Z$ has no positive jumps and, taking $\beta =1,$ solves Shi's conjecture in this particular case. It may seem surprising that the exponent $\theta$ shows no dependence on $\beta$, since the processes $A^{(\beta)}$ look different according as $\beta < 0$ or $\beta \ge 0.$ Notice however that this dependence is retrieved in the lower tail formulation: by self-similarity one has
$$\p\lcr \sup_{t\in [0, 1]} A^{(\beta)}_t < \eps  \rcr\;=\;
\eps^{(\alpha-1)/2(\alpha +\beta) +o(1)}, \qquad \eps \to 0,$$
as soon as $Z$ has no positive jumps.

Our proof of Theorem B is reminiscent of an argument of Sinai \cite{Si} who
had considered this question for $\beta =1,$ in discrete time. Roughly, we show
that $T$ belongs to an excursion interval of $Z$ whose left-end can
be controlled appropriately. In \cite{Si} - see the
end of Section 3 therein, the {\em simplicity} assumption on the
random walk is crucial to get the control. Here, we remark that this assumption can be
translated in continuous time into that of the absence of negative jumps, since then $X$ accumulates positive
value between $T$ and the next hitting time of $0$ by
$Z$. Alternatively, this accumulating argument is similar to that of Lemma
1 in \cite{Ber1}, which is key in analyzing the solution of the inviscid
Burgers equation whose initial data is a L\'evy process with no
positive jumps. The analogy is actually not very surprising, and we
will give more explanation for this at the end of the paper. 

To conclude this introduction, notice that Theorems A and B readily
entails the following result on the positive moments of $T$:

\begin{CORO} Suppose that $Z$ has no negative jumps. Then for every $k > 0,$
$$k < (\alpha-1)/2\alpha\quad\Longleftrightarrow\quad \e[ T^k]\; <\; +\infty.$$
\end{CORO}

In the third section of this paper, after proving Theorems A and B, we
will give some more results in this vein concerning the relevant
random variable $Z_T$. 

\section{Preliminaries}

\subsection{Recalls on stable processes and their homogeneous functionals} Consider $Z = \{Z_t, \; t\ge 0\}$ a real strictly $\alpha$-stable L\'evy process with index $\alpha\in (1,2]$, viz. a process with stationary and independent increments which is $(1/\alpha)-$self-similar:
$$\{ X_{kt}, \; t\ge 0\}\; \elaw\;\{ k^{1/\alpha}X_t, \; t\ge 0\}$$
for all $k >0$. Its L\'evy-Khintchine exponent $\Psi(\lambda) = - \log \e [\ee^{i\lambda Z_1}]$ is given by
\begin{equation}
\label{LK}
\Psi(\lambda) \; =\;\kappa\vert\lambda\vert^\alpha(1-{\rm i}\chi{\rm
  sgn}(\lambda)\tan(\pi\alpha/2)), \quad \lambda\in\r,
\end{equation}
where $\kappa > 0$ is the scaling parameter and $\chi\in [-1,1]$ is the skewness parameter. The positivity parameter $\rho\; =\; \p\lcr Z_t > 0\rcr$, independent of $t$ by self-similarity, is given by
$$\rho \; = \; \frac{1}{2}\; +\; (\pi\alpha)^{-1}\arctan (\chi\tan(\pi\alpha/2))$$
(see Theorem 2.6.3 in \cite{Zol}). Since $\alpha >1$, $\rho$ takes its values in the interval $[1-1/\alpha, 1/\alpha]$, the case $\rho = 1- 1/\alpha$ corresponding to the spectrally positive situation ($\chi = 1$ and $Z$ has no negative jumps) and the case $\rho = 1/\alpha$ to the spectrally negative situation ($\chi =-1$ and $Z$ has no positive jumps). Of course these values coincide when $\alpha = 2$ (the Brownian case).

Because ${\rm Re} (1/\Psi(\lambda) = O(\vert\lambda\vert^{-\alpha})$ as $\vert \lambda \vert\to \infty,$ Theorem V.1 in \cite{Ber} entails that $Z$ possesses a local time process $L=\lacc L(t,x), \; t\ge 0, \; x\in \r\racc$ in the sense that that for any non-negative Borel function $f$,
\begin{equation}\label{density}
\int_0^t f(X_s) \d s\; =\; \int_\r f(x) L(t,x) \d x, \qquad t\ge 0.
\end{equation}
Besides, we know from Exercise V.3 in \cite{Ber} that a.s. the map $x \mapsto L(t,x)$ is $\eta$-H\"older for every $\eta <(\alpha-1)/2$. This property and the formula (\ref{density}) allows to define the signed additive functionals
$$A^{(\beta)}_t\; =\; \int_0^t \vert Z_s\vert^\beta {\rm sgn} (Z_s) \, \d s, \quad t\ge 0$$ 
for every $\beta > - (\alpha +1)/2.$ For $\beta > -1,$ the integral
converges absolutely by the self-similarity of $Z$, because $\alpha
> 1$. For $\beta \in (-(\alpha +1)/2, 1],$ it is meant as a Cauchy principal value:
$$A^{(\beta)}_t\; =\;\lim_{\eps\to 0} \int_0^t \Un_{\{\vert Z_s\vert > \eps\}}\vert Z_s\vert^\beta {\rm sgn} (Z_s) \, \d s\; =\;\lim_{\eps\to 0}\int_\r \Un_{\{\vert x\vert > \eps\}} \vert x\vert^\beta {\rm sgn} (x) (L(t,x) -L(0,x)) \d x.$$
We refer to the seminal paper \cite{BY} for numerous results on the principal
values of Brownian local times, and also to \cite{FG} for their relations
with certain limit theorems in the general stable case.
 
\subsection{A symmetry lemma} Set $L_t = L(t,0)$ for simplicity and let $\lacc\tau_t, \; t\ge 0\racc$ be the right continuous inverse of $L_t$:
$$\tau_t \;= \;\inf\left\{ u\ge 0: \; L_u > t \right\} , \qquad t\ge 0.$$ 
From the strong Markov property for $Z$ at  times $\tau_t$, and the fact that
$$\tau_{t+s} - \tau_t \; =\; \inf\left\{ u\ge 0: \; L_{\tau_t + u} - t > s \right\}, \quad s, t \ge 0,$$
we see that the process $t\mapsto (\tau_t, A^{(\beta)}_{\tau_t})$ is a bivariate L\'evy process. Besides, it is an easy consequence of the self-similarity of $Z$ that $t\mapsto \tau_t$ is a stable subordinator with index $(\alpha-1)/\alpha$ and $t \mapsto A^{(\beta)}_{\tau_t}$ a stable process with index $(\alpha-1)/(\alpha +\beta).$ For the sake of concision, from now on we will use the notations $\gamma = (\alpha -1)/\alpha, \delta = (\alpha -1)/(\alpha +\beta)$, and set
$$\xi_t = A^{(\beta)}_{\tau_t}, \quad t\ge 0.$$ 
The following lemma, whose statement is trivial when $Z$ itself is symmetric, entails that $\xi$ is a symmetric $\delta$-stable process:

\begin{lemma}\label{sym} The random variables $(\tau_1, \xi_1)$ and $(\tau_1, -\xi_1)$ have the same law.
\end{lemma}

\noindent
{\it Proof.} Introducing the process
$$Z^\#_t \;=\;  \tau_1^{-1/\alpha}Z_{t\tau_1}, 
    \qquad t\in [0, 1],$$
we see after a change of variable that
$$\xi_1\; =\; \tau_1^{1+\beta/\alpha}\int_0^1 \vert Z^\#_s\vert^\beta {\rm sgn} (Z^\#_s) \, \d s.$$
On the other hand, it was established in \cite{DSS} that for any bounded measurable functional $F$,
 \begin{equation} \label{abslo-cty}
     \e\left[ F((Z^\#_t, \; t\in[0,1]), \tau_1) 
     \right]\; =\;  \e^{\rm (br)} \left[ 
     {c\over L_1} F\lpa (Z_t, \; t\in[0,1]), 
     L_1^{\alpha/(1-\alpha)} \rpa\right] 
 \end{equation}
where $\p^{\rm (br)}$ stands for the law of the standard bridge
associated with $Z$ (see section VIII.3 in \cite{Ber} for a
definition) and $c\in (0,\infty)$ is a normalization constant. In
passing, notice that (\ref{abslo-cty}) solves actually Exercise VIII.7
in \cite{Ber} and corrects a typographical error therein: the density of the law of the pseudo-bridge with respect to that of the standard bridge should be $kL(0,1)^{-1}$ and not $kL(0,1)^{1/(\alpha-1)}$. 

Setting $\theta_1 = \tau^{-1-\beta/\alpha}\xi_1$, we get from (\ref{abslo-cty}) that for every $f$ measurable bounded function 
$$\e\lcr f(\theta_1, \tau_1) \rcr\; =\;  \e^{\rm (br)} \lcr 
     {c\over L_1} f(A^{(\beta)}_1, L_1^{\alpha/(1-\alpha)})\rcr.$$
But from their very definition, 
$$(A^{(\beta)}_1, L_1^{\alpha/(1-\alpha)})\; =\; ({\bar A}^{(\beta)}_1, {\bar L}_1^{\alpha/(1-\alpha)})$$
where ${\bar A}^{(\beta)}$ and ${\bar L}$ are associated with the time-reversed process ${\bar Z} = \{ Z_{(1-t)-}, 0\le t\le 1\}$. From Exercise VIII.5 in \cite{Ber}, this entails that 
$$\e\lcr f(\theta_1, \tau_1) \rcr\; =\;  {\hat \e}^{\rm (br)} \lcr {c\over L_1} f(A^{(\beta)}_1, L_1^{\alpha/(1-\alpha)})\rcr$$
where ${\hat \p}^{\rm (br)}$ is the law of the standard bridge associated with the dual process ${\hat Z} = -Z.$ Using now (\ref{abslo-cty}) backwards, this yields 
$$\e\lcr f(\theta_1, \tau_1) \rcr\; =\; \e\lcr f(-\theta_1, \tau_1)\rcr$$
for every $f$ bounded measurable, so that $(\theta_1, \tau_1)$ and $(-\theta_1, \tau_1)$ have the same law, as well as $(\xi_1, \tau_1)$ and $(-\xi_1, \tau_1)$.

\qed

\begin{remarks} {\em (a) This lemma is intuitively true from the facts that the variable $(\tau_1, \xi_1)$ is invariant under the operation reversing the time of each excursion of $Z$ and making it right-continuous - see Definition (4.4) in \cite{GS}, and that the transformation of the law of $Z$ under this latter operation yields actually the law of ${\hat Z}$ - see Theorem (4.8) in \cite{GS}. Eventhough this argument is simpler, we feel that making it perfectly rigorous would force us to introduce some heavier notation than the one we use in the previous proof, and we leave this task to the interested reader. \

(b) It is possible to determine
the scaling parameter $\kappa_\tau = -\log \e\lcr
\ee^{-\tau_1} \rcr$ of the subordinator $\tau$ from Proposition V.4 in \cite{Ber}. More precisely, we have
$$\kappa_\tau \,= \, 2\pi\lpa \int_\r {\rm Re}\lpa
\frac{1}{1+\Psi(\lambda)}\rpa \d \lambda\rpa^{-1} \!\!\!\! = \; \pi\kappa^{1/\alpha}\lpa \int_0^\infty \lpa \frac{1 +\lambda^\alpha}{(1+\lambda^\alpha)^2 +\chi^2 \tan^2 (\pi\alpha/2)\lambda^{2\alpha}}\rpa \!\d \lambda\rpa^{-1}$$
where $\kappa$ is the scaling parameter of $Z$ appearing in (\ref{LK}). When
$\chi = 0$, using the 2nd formula p. 5 and the 3rd formula p. 10 in
\cite{Mag}, one obtains
$$\kappa_\tau \; =\; \alpha\kappa^{1/\alpha}\sin(\pi/\alpha).$$
When $\chi \neq 0$, using the 5th formula p. 188 and the 8th formula
p. 172 in \cite{Mag}, one gets
$$\kappa_\tau \; =\;
\frac{\alpha\kappa^{1/\alpha}\sin(\pi/\alpha)\mu^{-1/\alpha}\sqrt{1-\mu^2}}{\sin(\pi(\alpha-1)\vert\rho-1/2\vert)+\mu
  \sin(\pi\vert\rho-1/2\vert)}$$
where $\mu = (1+\chi^2\tan^2(\pi\alpha/2))^{-1/2}$ and $\rho$ is the
positivity parameter of $Z$. Notice that when $\vert\chi\vert =1,$ the formula becomes
$$\kappa_\tau \: =\;\frac{\alpha\kappa^{1/\alpha}}{(\sin((\alpha-1)\pi/2))^{1/\alpha}}\, ,$$
a fact which could have been obtained directly in solving Exercise VII.2 in \cite{Ber}.
 
However, it seems difficult to compute in general the scaling
parameter $\kappa_\xi = -\log \e\lcr
\ee^{{\rm i}\xi_1} \rcr$ of the symmetric stable process $\xi$.
In the Brownian case $\alpha = 2$, one can use the independence of the
positive and negative excursions and solve an homogeneous
Sturm-Liouville equation, which gives (see Formula (1 b) in \cite{BY})
\begin{equation}\label{constant}\kappa_\xi\; =\; \frac{\pi(2\kappa)^{1/\delta -1}2^\delta}{2\delta \sin (\pi\delta/2)}\lpa \frac{\delta^\delta}{\Gamma(\delta)}\rpa^2\cdot
\end{equation}
 When $\alpha < 2$ one cannot use such
a method because the excursions change a.s. their sign at least once during their lifetime. However, in the case $\beta = -1$, a remarkable formula due to Fitzsimmons, Getoor and Bertoin - see Theorem V.7 in \cite{Ber} - entails that
$$\e\lcr \ee^{{\rm i}\lambda A^{(-1)}_{\tau_t}} \rcr\; =\; \ee^{-t\pi\vert\lambda\vert},\quad \lambda\in\r$$
(notice that it is consistent with (\ref{constant}) above when $\delta =1$). The fact that $\kappa_\xi$ does not depend on $\kappa$ follows easily from the self-similarity of $Z,$ and the simple constant $\pi$ stems roughly from the identity ${\mathcal  F}{\mathcal  H} = \pi {\rm i}\, {\rm sgn}$ where ${\mathcal  F}$ (resp. ${\mathcal  H}$) is the Fourier (resp. Hilbert) transform. For $\beta = 1$, using the identity ${\mathcal  F}(x.) = -{\rm i}\, {\mathcal F}'$ and reasoning exactly as in section V.2 in \cite{Ber}, one can prove that the bivariate scaling parameter $\kappa (q, \lambda) = -\log \e [ \ee^{- q \tau_1 + {\rm i}\lambda \xi_1}]$ appears in the linear equation
$$\lambda\varphi'(t) \; +\; \kappa (q, \lambda)\; =\; (q +\Psi(-t))\varphi (t)$$
satisfied by the function $\varphi = \FF h,$ where $T_{-x} = \inf\{ t
>0, \; Z_t = -x\}$ 
and 
\begin{equation}
\label{chelles}
h(x) = \e [ \ee^{- q T_{-x} + {\rm i}\lambda A^{(1)}_{T_{-x}}}].
\end{equation}
Using Lemma VIII.13 in \cite{Ber}, one can show that $\vert x\vert
h(x) \in {\rm L}^1(\r)$, so that by (\ref{LK}) and the Riemann-Lebesgue lemma,
$$\kappa (q, \lambda) \; =\; \kappa(1 + \chi\tan (\pi\alpha/2))\lim_{t\to +\infty} t^\alpha \varphi(-t).$$ 
However this formula seems barely tractable, even in the Brownian case
where $h$ is explicit - see Formula (5') p. 387 in \cite{La}, but
given as a quotient of Airy functions. \

(c) When $\beta > -1$, the process $\xi$ has finite variations and hence, can be decomposed as the difference of two i.i.d. $\delta$-stable subordinators. There is actually another decomposition 
$$\xi\; = \;\xi^+ - \xi^-,$$
with the notation
$$\xi^+_t \; = \; \int_0^{\tau_t} \vert Z_s\vert^\beta \Un_{\{Z_s >0\}} \, \d s \quad{\rm and}\quad  \xi^-_t \; = \; \int_0^{\tau_t} \vert Z_s\vert^\beta \Un_{\{Z_s <0\}} \, \d s $$
for every $t\ge 0.$ Reasoning as above, one can show that the two
$\delta$-stable subordinators $\xi^+$ and $\xi^-$ have also the same law. Nevertheless they are not independent since they jump simultaneously, unless $\alpha = 2.$ }
\end{remarks}

\section{Proofs of the theorems}

\subsection{Proof of Theorem A} Let $\kappa(q,\lambda) = -\log \e [ \ee^{- q \tau_1 + {\rm i}\lambda \xi_1}]$ be the Laplace-Fourier exponent of the bivariate L\'evy process $\{(\tau_t, \xi_t), \, t\ge 0\}.$ Setting $I_t = \inf\{\xi_s, s\le t\}$ and $S_t = \sup\{\xi_s, s\le t\}$ for $t \ge 0$, the two-dimensional Wiener-Hopf factors
are defined followingly:
$$\phi_+(q,\lambda) \, =\, \kappa(q,0)\int_0^\infty \e\lcr \ee^{-q\tau_t + {\rm i}\lambda S_t}\rcr \d t \quad\mbox{and}\quad \phi_-(q,\lambda) \, =\, \kappa(q,0)\int_0^\infty \e\lcr \ee^{-q\tau_t + {\rm i}\lambda I_t}\rcr \d t$$
for all $q >0, \lambda\in \r.$ The Wiener-Hopf factorization due to
Isozaki - see Theorem 1 in \cite{Iso} - yields
$$\phi_+(q,\lambda)\phi_-(q,\lambda)\; =\; \frac{\kappa(q,0)}{\kappa(q,\lambda)}, \qquad q > 0, \lambda \in \r.$$
Notice that by Lemma \ref{sym}, the L\'evy processes $ (\tau_t,
\xi_t)$ and $(\tau_t, -\xi_t)$ have the same laws. This entails $\phi_-(q,\lambda) = {\bar \phi_+(q,\lambda)}$ for all $q, \lambda$, so that
\begin{equation}
\label{WH}
\vert\phi_+(q,\lambda)\vert\; =\; \sqrt{\frac{\kappa(q,0)}{\kappa(q,\lambda)}}, \qquad q > 0, \lambda \in \r.
\end{equation}
The factor $\phi_+$ is important for our problem. Indeed, setting 
$$\Theta (r) = \inf\{t > 0, \, \xi_t > r\},$$
it follows from the Markov property at time $\Theta (r)$ - see Lemma 1
and Lemma 2 in \cite{Iso} - that for every $r > 0$
$$1 - \e\lcr \ee^{-\tau_{\Theta(r)}}\rcr\; =\; \int_\r \lpa \frac{1 - \ee^{{\rm i }\eta r}}{2\pi {\rm i} \eta}\rpa \phi_+(1, \eta) \d \eta.$$
Since by scaling $\tau_{\Theta (r)} \elaw r^{\frac{\alpha}{\alpha +
    \beta}} \tau_{\Theta (1)}$, after a change of variable we obtain
\begin{equation}
\label{Laplace1}
1 - \e\lcr \ee^{-r\tau_{\Theta(1)}}\rcr\; =\; \int_\r \lpa \frac{1 - \ee^{{\rm i }\eta }}{2\pi {\rm i} \eta}\rpa \phi_+(1, \eta r^{-(1 + \beta/\alpha)}) \d \eta, \qquad r \ge 0.
\end{equation} 
We now proceed to an asymptotic analysis of $\phi_+(1, x)$ when $x\to\infty,$ mimicking Lemma 3 in \cite{Iso}. First, by (\ref{WH}), a scaling argument and the symmetry of $\xi$, \begin{equation}
\label{module}
\vert\phi_+(1,x)\vert\; =\; \sqrt{\frac{\kappa(1,0)}{\kappa(1,x)}}\; =\; \sqrt{\frac{\kappa_\tau\vert x\vert^{-\delta}}{ \kappa(\vert x\vert^{-\alpha/(\alpha + \beta)}, 1)}}\; \sim\;\sqrt{\frac{\kappa_\tau}{\kappa_\xi}}\vert x\vert^{-\delta/2}, \quad x\to\infty.
\end{equation}
Second, by Theorem 1 in \cite{Iso} - noticing that $(\ee^{{\rm i }\xi_t -1})$ should be read $(\ee^{{\rm i }\xi_t} -1)$ therein, a scaling argument and a change of variable,
\begin{eqnarray}
\label{argument}
\arg (\phi_+(1,x)) & = & \int_0^\infty \frac{\d t}{t} \e\lcr \ee^{-\tau_t} \sin (x \xi_t) \Un_{\{ \xi_t > 0\}}\rcr\nonumber\\
& = & {\rm sgn} (x) \int_0^\infty \frac{\d t}{t} \e\lcr \ee^{-t^{1/\gamma}\tau_1} \sin (\vert x\vert t^{1/\delta}\xi_1) \Un_{\{ \xi_1 > 0\}}\rcr\nonumber\\
& = & \delta \,{\rm sgn} (x) \int_0^\infty \frac{\d u}{u} \e\lcr \ee^{-\tau_1\lpa \frac{u}{\vert x\vert}\rpa^{\alpha/(\alpha +\beta)}} \sin (\vert u\vert \xi_1) \Un_{\{ \xi_1 > 0\}}\rcr,\nonumber\\
& \sim & \delta \,{\rm sgn} (x) \int_0^\infty \frac{\d u}{u}  \e\lcr\sin (\vert u\vert \xi_1) \Un_{\{ \xi_1 > 0\}}\rcr\; = \; \frac{\delta\pi}{4} \,{\rm sgn} (x),\quad x \to\infty.
\end{eqnarray}
Last, we compute an oscillating integral which will appear at the
limit. Notice that this result is also quoted in \cite{Iso} p. 223,
nevertheless we will include a proof for the sake of clarity.

\begin{lemma} \label{oscil}
One has
$$\int_\r \lpa \frac{1 - \ee^{{\rm i }\eta}}{2\pi {\rm i} \eta}
\rpa\vert\eta \vert^{-\delta/2}\ee^{{\rm i }\delta {\rm sgn }(\eta)\pi/4} \d \eta \; =\; \frac{2}{\delta \Gamma(\delta/2)}\, \cdot$$
\end{lemma}
\noindent
{\em Proof.} Setting $I$ for the value of the integral we can rewrite,
by a parity argument,
$$I\; =\; \int_\r \lpa \frac{\sin (\eta/2)}{\pi  \eta}
\rpa\vert\eta \vert^{-\delta/2}\ee^{{\rm i }(\delta {\rm sgn }(\eta)\pi/4-\eta/2)} \d \eta\; =\;  \frac{2}{\pi}\int_0^\infty \sin (\eta/2)\eta^{-1-\delta/2}\cos(\delta\pi/4-\eta/2)\d \eta.$$
After trigonometric transformations and integrating by parts, we get
\begin{eqnarray*}
I & = & \frac{\cos (\delta\pi/4)}{\pi}\int_0^\infty\eta^{-1-\delta/2}\sin \eta \d \eta \; +\; \frac{2\sin (\delta\pi/4)}{\pi\delta}\int_0^\infty\eta^{-\delta/2}\sin \eta\d \eta \\
& = & \sin(\delta\pi/2)\lpa \frac{\Gamma(1-\delta/2)}{\pi\delta} - \frac{\Gamma(-\delta/2)}{2\pi}\rpa \; =\; \frac{2}{\delta\Gamma(\delta/2)}
\end{eqnarray*}
where we used the 5th formula p. 9 in \cite{Mag} in the second
equality, and the formulae $\Gamma(z)\Gamma(-z) = -\pi/z\sin(\pi z), \Gamma(z)\Gamma(1-z) = \pi/\sin(\pi z)$ in the third.

\qed

\vspace{2mm}

Notice by the assumption on $\beta$, one has $r^{-(1+\beta/\alpha)}\to
+\infty$ as $r\to 0.$ Hence, putting (\ref{Laplace1}), (\ref{module}),
(\ref{argument}) and Lemma \ref{oscil} together one obtains
$$1 - \e\lcr \ee^{-r\tau_{\Theta(1)}}\rcr\;
\sim\;\sqrt{\frac{\kappa_\tau}{\kappa_\xi}} \lpa\frac{2}{\delta\Gamma(\delta/2)} \rpa r^{\gamma/2}, \qquad r\to
0,$$
the dominated convergence argument being plainly justified by
(\ref{Laplace1}) and (\ref{module}). By a Tauberian theorem and a monotone density theorem - see
e.g. \cite{Ber} p. 10, this yields 
\begin{equation}
\label{tau}
\p[\tau_{\Theta(1)} > t]\;\sim\; \sqrt{\frac{\kappa_\tau}{\kappa_\xi}}
\lpa\frac{2}{\delta\Gamma(\delta/2)\Gamma(1-\gamma/2)}\rpa t^{-\gamma/2},
\qquad t\to +\infty.
\end{equation}
The key-point is now that by the definition of $T$, for every $t > 0$ 
$$\{T > t\}\; =\; \{A^{(\beta)}_s < 1, \; s \le t\}\;\subset \;\{\xi_u < 1,
\,\forall\, u \; / \; \tau_u \le t\}\; =\; \{ \Theta_1 > L_t\}\;\subset\{\tau_{\Theta(1)} > t\}.$$
This entails finally
$$\p[T > t]\; \le \; \k\, t^{-\gamma/2}, \quad t\to +\infty,$$
with 
$$  \k\; =\; \sqrt{\frac{\kappa_\tau}{\kappa_\xi}}
\lpa\frac{2}{\delta\Gamma(\delta/2)\Gamma(1-\gamma/2)}\rpa$$
explicit when $\alpha =2$ or $\beta = -1,$ thanks to the above
Remark 2 (a).

\qed

\begin{remark}{\em At a less precise level and for the case 
where $\beta =1$ and $Z$ is symmetric, it had been shown in \cite{DSS} that 
$$\p\lcr \sup_{t\in [0, 1]} A^{(1)}_t < \eps  \rcr\;\le\;
\eps^{(\alpha-1)/2(\alpha +1) +o(1)}, \qquad \eps \to 0$$
with a different method: first, one remarks that
$$\p\lcr \sup_{t\in [0, \tau_1]} A^{(1)}_t < \eps  \rcr\;\le\;\p\lcr \sup_{t\in [0, 1]} A^{(1)}_{\tau_t} < \eps  \rcr\;\sim\;
k \eps^{(\alpha-1)/2(\alpha +1)}, \qquad \eps \to 0$$
for some finite positive constant $k,$ because $t\mapsto A^{(1)}_{\tau_t}$ is a symmetric $(\alpha-1)/(\alpha +1)$-stable process and in view of Proposition VIII.3 in \cite{Ber}. Second, at the cost of some $o(1)$ term on the exponent, one replaces $\tau_1$ by $1$ in the probability on the left-hand side after a suitable slicing of the big values of $\tau_1.$}
\end{remark}
\subsection{Proof of Theorem B}

Set $\Theta = \inf\{ t > 0, \, \xi_t \ge 1\}$ and for every $c > 0$
introduce the function
$$f^c(t)\; =\; (c \log t)^{1/(1-\alpha)}, \qquad t >1.$$ 
One has
\begin{eqnarray*}
\p\lcr \Theta > t^\gamma\rcr & \le & \p\lcr \tau_{t^\gamma/2} < t
f^c(t)\rcr\; +\; \p\lcr \Theta > t^\gamma, \, \tau_{t^\gamma/2} \ge t
f^c(t)\rcr\\
& \le & \p\lcr \tau_{t^\gamma/2} < t
f^c(t)\rcr\; +\; \p\lcr \tau_{\Theta-} \ge t
f^c(t)\rcr
\end{eqnarray*}
since $\tau$ is a.s. increasing. The key-point is the inequality 
\begin{equation}
\label{KP}
T \; \ge\; \tau_{\Theta-}\quad\mbox{a.s.}
\end{equation}
Indeed, if we set $T' = \inf\{t \ge T, \, Z_t = 0\},$ then it follows
from the absence of negative jumps for $Z$ that $Z_T \ge 0$ and
$A^{(\beta)}_{T'} \ge 1$ a.s. Hence, since a.s. $T\le\tau_\Theta$ we
see that a.s. $T' = \tau_\Theta$, which readily entails (\ref{KP}) by
the definition of $T'$. We deduce
\begin{eqnarray*}
\p\lcr T \ge t f^c(t)\rcr & \ge & \p\lcr \Theta > t^\gamma\rcr \; -\;
\p\lcr \tau_{t^\gamma/2} < t f^c(t)\rcr\\
& \ge & \p\lcr \Theta > t^\gamma\rcr \; -\;
\p\lcr \tau_{1/2} < f^c(t)\rcr\\
& \ge & k_1 t^{-\gamma/2}\; -\; t^{-ck_2},\qquad t \to \infty
\end{eqnarray*}
for some $k_1, k_2 >0,$ where in the third inequality
we used Proposition VIII. 2 and the remark after
Proposition VIII.3 p. 221 in \cite{Ber} - recall that $\tau$ is a $\gamma$-stable subordinator. Taking $c$ sufficiently large, we
obtain
$$\p\lcr T \ge t f^c(t)\rcr\; \ge\; k_3 t^{-\gamma/2},\qquad t
\to\infty$$
for some positive constant $k_3$. Changing the variable $u = tf^c(t)$,
we see by the definition of $f^c$ that $t \le u(\log
u)^{1/(\alpha-1)}/2$ provided $u$ is sufficiently large. Putting
everything together yields
$$\p\lcr T \ge u \rcr\; \ge\; k_4 u^{-\gamma/2}(\log u)^{-1/2\alpha},\qquad u
\to\infty$$
for some constant $k_4 >0,$ as required.

\qed

\begin{remarks}{\em 

(a) In the above proof, the key-argument (\ref{KP}) is completely peculiar to the spectrally positive framework. When there are negative jumps, it seems possible that 
$$T' \; \ll\; \tau_{\Theta}$$
with significant probability, and this raises the question whether the
exponent $\gamma/2$ is still the critical one. In a recent paper \cite{TS} we actually proved that it is not the case anymore when there are negative jumps, at least when $\alpha$ is close to 1, hence contradicting the validity of Shi's conjecture in general. At the end of this paper, we will state another conjecture on the value of the critical exponent $\theta$ for integrated L\'evy $\alpha$-stable processes.\

(b) In the Brownian case $\alpha = 2$, it is possible to obtain a
sharper bound
$$\p\lcr T > t \rcr\; \ge \; \k' t^{-1/4}, \qquad t\to\infty$$
for some explicit constant $\k'$ - see Theorem 3 in \cite{Iso}. The
method consists in proving that
$$\p\lcr \tau_\Theta - T > t\rcr \;\le\;\ \k''t^{-1/4},\qquad t\to\infty,$$
where $\k''$ is made explicit with the help of the function $h$ in
(\ref{chelles}), and shown to be strictly smaller than
$\k$. Notice that the arguments of \cite{Iso} are given only for 
$\beta > 0.$ Nevertheless they can be extended to the general situation $\beta > -3/2$ with the same formula for $\k$ and $\k'$. As mentioned in the introduction, in a subsequent work \cite{IK} it is proved  that 
$$\p\lcr T> t\rcr\; \sim\; \k_\beta t^{-1/4}, \qquad t\to \infty$$
as soon as$\beta > -1,$ for some explicit constant $\k_\beta$ - see Corollary 1 in \cite{IK}. Notice that the value $\k_1 =\; 3^{4/3}\Gamma(2/3)/\pi 2^{13/12}\Gamma(3/4)$ had been previously computed by Goldman \cite{Go}. \

(c) In the Non-Gaussian case, it seems
difficult to obtain sharper estimates with the above method, due to the
lack of information on $h$. When $\beta \ge 0$ however, it is possible to
obtain an explicit uniform bound 
\begin{equation}
\label{uniform}
\p\lcr \tau_\Theta - T > t \rcr\; \le \; \k\, t^{-\gamma/2}, \qquad
t \ge 0
\end{equation}
thanks to the formula
\begin{equation}
\label{skoro}
\p\lcr \tau_\Theta - T > t \rcr \; = \; \alpha\!\int_0^\infty\!\! {\widehat
  p}_1(u)\,\p [ Z_T > t^{1/\alpha}u ] \d u
\end{equation}
where ${\widehat p}_1$ is the density of ${\widehat Z}_1 = -Z_1,$ and an upper bound on
 $\p\lcr Z_T > u\rcr$ which will be detailed soon afterwards. The
formula (\ref{skoro}) follows at once from the Markov property at time $T$, an
integration by parts, a scaling argument and the well-known fact that
\begin{equation}
\label{bingham}
\sup\{{\widehat Z}_t, \, t\le 1\} \;\elaw \;{\widehat Z}_1\Un_{\{
  {\widehat Z}_1 \ge 0\}}
\end{equation}
since ${\widehat Z}$ has no positive jumps - see e.g. \cite{Bi}
p. 749. Unfortunately, it seems that (\ref{uniform}) is useless as soon as we do not have any more precise information on $\kappa_\xi$. \

(d) Recently, the simple method of Theorem B has been applied successfully to derive some a priori difficult lower bound on the small ball probabilities for a class of homogeneous Brownian functionals \cite{LaSi}.
}
\end{remarks}

\subsection{Final remarks and conjectures} In this paragraph we first give some comments on the positive moments of the interesting random variable $Z_T$. In the Brownian case and for $\beta =1,$ the law of $Z_T$ had been computed by Gor'kov - see Formula (3) in \cite{La}, but the expression seems too complicated to allow further computations. In the general case with no negative jumps, using (\ref{skoro}), the formula 
$$\e[S^k]\; =\; k\int_0^\infty s^{k-1} \p [S > s] \d s$$
which is valid for any positive random variable $S$ and any $k >0,$ and Fubini's theorem, one obtains 
\begin{eqnarray*}
\e[(\tau_\Theta - T)^k] & = & \alpha\lpa\int_0^\infty u^{-k\alpha } {\widehat p}_1(u) \d u\rpa \e[Z_T^{k\alpha }]\\
& = & \lpa \frac{\Gamma (1-k\alpha)}{\kappa^k\Gamma(1-k)}\rpa\e[Z_T^{k\alpha }]
\end{eqnarray*}
for every $k \in [0, 1/\alpha),$ where the computation of the integral in the second line comes from Theorem 2.6.3 in \cite{Zol}. In the Brownian case, using (5.6) in \cite{Iso} yields then the following criterion:
$$ k < 1/2\quad\Longleftrightarrow\quad\e[ Z_T^k]\; <\; +\infty,$$
for every $k > 0.$ In the general stable case with no negative jumps, using Theorem A we obtain the inclusion 
\begin{equation}
\label{inclu}
k < (\alpha-1)/2\quad\Longrightarrow\quad \e[ Z_T^k]\; <\; +\infty,
\end{equation}
which we actually believe to be an equivalence. Notice that when $\beta\ge 0,$ one can also obtain the following reinforcement of (\ref{inclu}):
$$\p[Z_T > u] \; \le\; \k\, u^{-(\alpha -1)/2}, \quad u\to \infty,$$
for some explicit constant $\k$, which together with Theorem 2.6.3 in \cite{Zol} leads to the aforementioned estimate (\ref{uniform}). Let me sketch the proof. First, by the Markov property at time $T$, a scaling argument and Skorokhod's formula (\ref{bingham}) one gets
$$\p[Z_T > u] \; =\; \frac{\p[Z_T > u, \, T_{cZ_T}\circ\theta_T > Z_T^\alpha]}{\alpha \p[0< {\widehat Z}_1 < 1-c]}$$
for every $c \in (0,1)$ and $u >0,$ where $\theta$ is the shift
operator of $Z$ and $T_x$ the first hitting time of $x$ by $Z$. Second, since $Z$ has no negative jumps and since $\beta \ge 0,$ one deduces
$$\p[Z_T > u] \; \le\; \frac{\p[A^{(\beta)}_{\Theta(1)} > c u^{\alpha +\beta}]}{\alpha \p[0< {\widehat Z}_1 < 1-c]}$$
for every $c \in (0,1)$ and $u >0.$ Last, remembering that $A^{(\beta)}$ is a symmetric $\delta$-stable process and using Exercise VIII. 3 in \cite{Ber}, one obtains
$$\p[Z_T > u] \; \le\; \k_c \, u^{-(\alpha-1)/2}, \quad\, u\to\infty$$
for an explicit constant $\k_c$ which can be minimized in $c\in(0,1).$ \\

In the general case with or without negative jumps, finding a lower bound on the quantities $\p[Z_T > u]$ seems at least as difficult as obtaining a lower bound on $\p[T > u]$. In the case without {\em positive} jumps, the two asymptotics are clearly connected thanks to the basic inequality
$$\p[Z_T > u]\; \le\; \p[T > u^{\alpha(1 - o(1))}]\; +\; \p[S_1 >
u^{o(1)}],\quad u\to\infty,$$
and the fact that $S_1 =\sup\{ Z_t, \, t\le 1\}$ has finite exponential moments - see Corollary VII.2 in \cite{Ber}. From the L\'evy-Khintchine formula and the optional sampling theorem at time $T,$ one can finally raise the following

\begin{CONJ} For every $k >0$ one has
$$\e[ T^k]\; <\; +\infty\quad\Longleftrightarrow\quad\e[ Z_T^{k\alpha}]\; <\; +\infty.$$
\end{CONJ}
To put it in a nutshell, in this paper we have showed that the critical value of $k$ in the above conjecture is $(\alpha -1)/2\alpha$ when there are no negative jumps, but we still do not know what this value should be when there are negative jumps.
It might be worth mentioning that if we set $S = \inf\{ t > 0, \; Z_t > 1 \},$ then from Proposition VIII. 2 and Exercise VIII. 3 in \cite{Ber} we have
$$ k < \rho\quad\Longleftrightarrow\quad\e[ S^k]\; <\; +\infty\quad\Longleftrightarrow\quad\e[ Z_S^{k\alpha}]\; <\; +\infty$$
as soon as $Z$ has positive jumps (when $Z$ has no positive jumps only the first equivalence remains valid since $Z_S = 1$ a.s.) From this fact and the results of the present article, one might suppose that in general the critical exponent is
\begin{equation}\label{k0}
k_0\; =\; \frac{\rho}{2}
\end{equation}
in the above conjecture - recall that $\rho = (\alpha -1)/\alpha$  when $Z$ has no negative jumps and $\alpha >1$. Besides, when $\beta =1$, one can show from Proposition 3.4.1 in \cite{ST} that 
$$\p[ A^{(1)}_1 > 0] =\rho,$$
and it would be somewhat surprising that this quantity does not play any r\^ole in the lower tails of $A^{(1)},$ since it does in the lower tails of $Z$. \

Let us finally mention that if (\ref{k0}) held, then it would also be very interesting to know if for $\beta = 1$ the following upper bound dealing with the {\em non-strictly} stable process $Z' = \{Z_t-t, \, t\ge 0\},$ remains true: 
\begin{equation}
\label{Burgers}
\p\lcr \int_0^t Z_s \d s < \eps \; + \; t^2,\;\forall\, t\in [0,1] \rcr\;\le\;\eps^{\rho/2 +o(1)}, \qquad \eps \to 0.
\end{equation}
Indeed, following a seminal approach due to Sinai \cite{Si1}, Molchan and Kholkhov recently showed that upper estimates like (\ref{Burgers}) provide the key-step to get upper bounds on the Hausdorff dimension of the so-called regular Lagrangian points for the inviscid Burgers equation with random initial data - see Theorem 1 in \cite{MK}. Molchan and Kholkhov's arguments are given in a fractional Brownian framework, but one can show \cite{TS} that
they remain valid for L\'evy stable processes: if (\ref{Burgers})
held, then one would have 
$${\rm Dim}_H \,\LL\; \le\; 1 - \rho\; \le\; 1/\alpha\quad\mbox{a.s.}$$
where $\LL$ is the set of Lagrangian regular points for the inviscid
Burgers equation with initial data $\widehat{Z}$. We refer to \cite{Si1, Ber1, JW} for various presentations of the Burgers
equation with random initial impulse and we recall that in \cite{JW}, the following well-known conjecture is made:
\begin{equation}
\label{JaWo}
{\rm Dim}_H \, \LL\; = \; 1/\alpha\quad\mbox{a.s.}
\end{equation}
This conjecture is true for Brownian motion by the result of Sinai \cite{Si1} - see also \cite{Ber1} for more rigorous arguments and more general results, and when ${\widehat Z}$ has no positive jumps i.e. $Z$ has no negative jumps and $\rho = 1 -1/\alpha$ - see again \cite{Ber1}. The situation where ${\widehat Z}$ has positive jumps seems to require a completely different methodology - see the conclusion of \cite{Ber1}, which we believe to be  connected with the optimal upper bound on $\p[T>u]$ when $Z$ has negative jumps: in a related paper \cite{TS}, we proved that in this situation (\ref{JaWo}) is false at least when $\alpha$ is close to 1, and we conjectured that
$${\rm Dim}_H \, \LL\; = \; 1 -\rho\quad\mbox{a.s.}$$
with the above notations, which is clearly connected to (\ref{k0}) as far as the upper bound is concerned.

\bigskip 
\noindent {\bf Acknowledgements.} I am grateful to A.~Devulder and Z.~Shi for leaving me the hands free to work alone on this paper. Thanks also go to J. Bertoin for providing me the reference \cite{GS}.

\end{document}